\documentclass[11pt]{amsart}
\usepackage{amsfonts,amssymb,amsmath}
\usepackage[all]{xy}
\begin{document}

\newcommand{\half}{{\textstyle{\frac{1}{2}}}}
\newcommand{\quarter}{{\textstyle{\frac{1}{4}}}}
\newcommand{\C}{{\mathbb C}}
\newcommand{\G}{{\mathbb G}}
\newcommand{\Gl}{{\rm Gl}}
\newcommand{\bG}{{\bf G}}
\newcommand{\bgamma}{{\boldsymbol{\Gamma}}}
\newcommand{\bH}{{\mathbb H}}
\newcommand{\BH}{{\bf H}}
\newcommand{\Q}{{\mathbb Q}}
\newcommand{\R}{{\mathbb R}}
\newcommand{\Z}{{\mathbb Z}}
\newcommand{\bz}{{\bf z}}
\newcommand{\odd}{{\rm odd}}
\newcommand{\ev}{{\rm even}}
\newcommand{\res}{{\rm res}}
\newcommand{\Sp}{{\rm Sp}}
\newcommand{\ess}{{\sf s}}
\newcommand{\SO}{{\rm SO}}
\newcommand{\Oh}{{\rm O}}
\newcommand{\T}{{\mathbb T}}
\newcommand{\U}{{\rm U}}
\newcommand{\KO}{{\rm KO}}
\newcommand{\rH}{{\rm H}}

\title {An integral lift of the $\Gamma$-genus}
\author{Jack Morava}
\address{Department of Mathematics, Johns Hopkins University, Baltimore,
Maryland 21218}
\email{jack@math.jhu.edu}
\thanks{The author was supported in part by the NSF}
\date {15 July 2012}

\begin{abstract}{The Hirzebruch genus of complex-oriented manifolds associated to
Euler's $\Gamma$-function lifts to a homomorphism of ring-spectra associated to
a family of deformations of the Dirac operator, parametrized by the homogeneous space 
$\Sp/\U$.}\end{abstract}

\maketitle

\noindent{\bf Introduction} \bigskip
                 
\noindent
Kontsevich, in his early work on deformation quantization [12 \S 4.6], drew attention to interesting formal
properties of Euler's $\Gamma$-function, regarded as defining something like a Hirzebruch genus.
This note presents that idea in the language of cobordism and formal groups, following [16]. The formalism 
of multiplicative power series defines a homomorphism
\[
\chi_\infty \; : \; M\U_* \to \C[v]
\]
(of graded rings, with a book-keeping indeterminate $v$) having no very immediate integrality properties, 
but classical function theory [\S 2.3.1] shows it to take values in the ring $\Q[\tilde{\zeta}(\odd)]$ generated
over the rationals by normalized zeta-values, usually expected to be transcendental. The principal 
result here [\S 3.1] is that a topologically reasonable homomorphism
\[
\xymatrix{
M\U \ar[r]^<<<<{\bgamma} & M\U \wedge_{M\Sp} \KO \ar[r]^{\cong \; [\half]} & \Sp/\U \wedge \KO [\half]
} 
\]
of ring-spectra provides a lift of $\chi_\infty$, via the composition
\[
\xymatrix{
\KO_*(\Sp/\U) \ar[r]^<<<<<{ch} & H_*(\Sp/\U,\Q[v^{\pm 1}]) \ar[r] & H_*(B\U,\Q[v^{\pm 1}]) \ar[r] & \C[v^{\pm 1}] }
\]
which sends primitive generators of $H_*(\Sp/\U,\Q)$ to odd $\zeta$-values. \bigskip

\noindent
It was the appearance of these periods (and their relation to the theory of mixed Tate motives in algebraic 
geometry) that precipitated much of the interest in the $\Gamma$-genus. They appear in the lift as 
generic parameters for a family of deformations of a Dirac operator over the homogeneous space $\Sp/\U$.
This seems to have interesting connections with [10] and [17]. \bigskip

\noindent
I'd like to thank Professor Hirzebruch for interest and conversation about this material, and Peter Landweber
and Ulrike Tillmann for helpful correspondence; but I owe special thanks to Bob Stong, for watching over my shoulder 
as I wrote. \bigskip

\section{Coigns of vantage}
\bigskip

\noindent
{\bf 1.0} It's useful to distinguish a coordinate $z$ at a point $x_0$ of a space $X$ from the
corresponding parametrization of a neighborhood $U \ni x_0$ : the former is a nice function
\[
\xymatrix{
{X \supset U}  \ar[r]^z & A }
\]
sending $x_0$ to 0 in some commutative ring $A$, while the latter is the map
\[
\bz : {\rm Spec} \; A \to U \subset X
\]
it defines (assuming we're in a context where this makes sense). \bigskip

\noindent
{\bf 1.1} For example, at the point $x_0 = [1:1]$ of the projective line, we have a coordinate 
\[
[u:1] \mapsto u - 1 := z 
\]
which defines the parametrization
\[
z \mapsto [1+z:1] 
\]
of a neighborhood of $[1:1]$. Similarly, 
\[
[q:1] \mapsto q^{-1} := z 
\]
is a coordinate at $[1:0] = \infty \in P_1$, while
\[
[x:1] \mapsto x := z 
\]
is a coordinate at $[0:1] = 0$. \bigskip

\noindent
{\bf 1.2} An abelian group germ $\G$ at $x_0 \in X$ is the germ of a function
\[
\G : U \times U, x_0 \times x_0 \to U,x_0
\]
satisfying identities such as 
\[
\G(x,\G(y,z)) = \G(\G(x,y),z), \; \G(x,x_0) = \G(x_0,x) = x, \; \&c \;;
\]
if $\G$ is suitably analytic, then a coordinate $z$ at $x_0$ associates to $\G$, the 
formal group law
\[
(z \circ \G)(\bz \times \bz) :=  z_0 +_\G z_1 \in A[[z_0,z_1]] \;.
\]
For example, the additive group germ $\G_a(x,y) = x+y$ at $[0:1] \in P_1$ defines $z_0,z_1
\mapsto z_0 + z_1$, while the multiplicative group germ $\G_m(u,v) = uv$ at $[1:1]$ defines
\[
z_0 +_{\G_m} z_1 = z_0 + z_1 + z_0z_1 
\]
(with coordinates as above). Different choices of coordinate (for fixed $\G$ and $x_0$) define, 
in general, distinct (but isomorphic) formal group laws: for example, if $t \in A^\times$ then 
$z = t^{-1}(u - 1)$ associates the formal group law
\[
z_0,z_1 \mapsto z_0 + z_1 + tz_0z_1\;.
\]
to the multiplicative group at $[1:1]$. \bigskip

\noindent
{\bf 1.3.1} The introduction of such a variable $t$ suggests the consideration of families, or
deformations, of group laws:
\[
u,v \mapsto \frac{uv}{1 - t(u-1)(v-1)}
\]
at $[1:1]$ (easily checked, eg for nilpotent $t$, to satisfy the axioms) is an interesting
example. With coordinate as above, the associated group law 
\[
z_0,z_1 \mapsto \frac{z_0 + z_1 + (1+t)z_0z_1}{1 - tz_0z_1} \;;
\]
is (strictly) isomorphic to $+_{\G_m}$, under the coordinate change
\[
z \to (1+t)^{-1} \log \left[\begin{array}{cc}
                                                     t & 1 \\
                                                    -1 & 1
                                               \end{array}\right](z) \in \Q[t][[z]] \;;
\]
note that the fractional linear transformation fixes $[1:1]$. \bigskip

\noindent
{\bf 1.3.2} Similarly, $\exp_A(z) :=  2\sinh z/2$ defines 
\[
z_0 +_A z_1 = z_0(1 + \quarter z_1^2)^{1/2} + z_1(1 + \quarter z_0^2)^{1/2} \in \Z[\half][[z_0,z_1]] \;,
\]
which is a specialization (at $\delta = - \textstyle{\frac{1}{8}}, \; \epsilon = 0$) of the formal group
law
\[
z_0 +_E z_1 = \frac{z_0 R(z_1) + z_1 R(z_0)}{1 - \epsilon z_0^2 z_1^2}
\]
defined by Jacobi's quartic $Y^2 = R(X)^2 := 1 - 2\delta X^2 + \epsilon X^4$.  
\bigskip

\noindent
{\bf 1.4} The focus of this note is the group germ
\[
\G_\infty : [q_0:1], \; [q_1:1] \mapsto [\Gamma(\log_\infty(q_0^{-1}) + \log_\infty(q_1^{-1})):1]
\]
at $\infty \in P_1(\R)$ defined by the expansion
\[
\exp_\infty(z) := z \; \exp(\gamma z - \sum_{k \geq 2} \frac{\zeta(k)}{k} (-z)^k) \in \R[[z]]
\]
of the entire function $\Gamma(z)^{-1}$ near 0 (with $\log_\infty(z)$ denoting its formal 
composition inverse): thus 
\[
z_0 +_{\G_\infty} z_1 = \Gamma(\log_\infty(z_0) + \log_\infty(z_1))^{-1} = z_0 + z_1 + 
2\gamma z_0z_1 + \dots \in \R[[z_0,z_1]] 
\]
with $z_k = q_k^{-1}$. Ohm's law for parallel resistors\begin{footnote}{a.k.a. the harmonic
mean of Archytas of Tarentum}\end{footnote}, in comparison, defines a group germ
\[
[q_0:1],[q_1:1] \mapsto [1:q_0^{-1} + q_1^{-1}]
\] 
at $\infty$, which (because $\frac{xy}{x+y}$ is not differentiable at $(0,0)$) 
is not analytic. \bigskip

\section{Characteristic classes and Hirzebruch genera} \bigskip 

\noindent
{\bf 2.1} A complex line bundle $\lambda \in H^1(X,\C^\times)$ has an associated
class
\[
\lambda^{-1} d \lambda \mapsto 2 \pi i [\lambda] : H^1(X,\Z(1)) \to H^2(X,2\pi i\Z)
\]
corresponding to the coordinate [1 \S 2.3, 19 \S 5.10]
\[
z = vx \in H^\ev(X,\Z[v^{\pm 1}])
\]
on the Picard group of topological complex line bundles. Interpreting $v$ as the product of 
the Bott class with Deligne's motive $2 \pi i$ reconciles some conventions of algebraic 
geometry with those of algebraic topology: for example 
\[
\frac{\pi [\lambda]}{\sin \pi [\lambda]} \mapsto \frac{vx/2}{\sinh vx/2} \;.
\]
When the grading is of background interest, I'll set $v$ equal to 1. \bigskip

\noindent
{\bf 2.2.1} A (one-dimensional) formal group law over a $\Q$-algebra $A$ can be written uniquely
as
\[
z_0 +_\G z_1 \; = \; \exp_\G (\log_\G(z_0) + \log_\G(z_1)) \;;
\]
in that case let
\[
H_\G(z) \; := \; \frac{z}{\exp_\G(z)} \in A[[z]]^\times
\]
denote its Hirzebruch multiplicative series [8 \S 15.5]. The function
\[
M \mapsto \Big(\prod_{i=1}^{i=n} H_\G(vx_i)\Big)[M] \; \in A[v]
\]
from (cobordism classes of) compact closed complex-oriented manifolds
of real dimension $2n$, with Chern roots $x_i$, defines a homomorphism
\[
\chi_\G : M\U_* \to A[v]
\]
of graded rings: the Hirzebruch genus associated to the group law $\G$. By a theorem 
of Mishchenko,
\[
\log_\G(v) = \sum_{n \geq 1} \frac{\chi_\G(P_{n-1}(\C))}{n} \in A[[v]] \;;
\]
the deformation of the multiplicative group in \S 1.3.1, for instance, represents 
Hirzebruch's genus $\chi_{-t}$ genus (defined on smooth projective complex varieties by 
\[
V \mapsto \sum (-1)^p (-t)^q \; \dim_\C H^{p,q}_{\rm dg}(V) \; v^{\dim_\C V} \;).
\] 
The coordinate rescaling $v \mapsto t^{-1/2}v$ sends its logarithm to
\[
\sum_{n \geq 1} [n](t)\frac{v^n}{n} 
\]
(with Gaussian $\frac{t^{n/2} - t^{-n/2}}{t^{1/2} - t^{-1/2}} = [n](t)$), and its formal group law to
\[
X,Y \mapsto \frac{X + Y + (t^{1/2} + t^{-1/2})vXY}{1 - vXY}
\]
(which is symmetric under the involution $t \mapsto 1/t$). \bigskip

\noindent
{\bf 2.2.2} I'll refer below to $M\SO, \; M\U$, and $M\Sp$ as the cobordism theories 
of $\R, \; \C$, and $\bH$-oriented manifolds, respectively. \bigskip

\noindent
The Pontryagin classes 
\[
p^\SO_t(V) = \sum_{k \geq 0} p^\SO_k(V) t^{2k} := \sum_{k \geq 0} (-1)^k c_{2k}(V \otimes \C) t^{2k}
\]
of a real vector bundle $V$ are defined in terms of the Chern classes of its complexification;
if $V$ was complex to begin with, then 
\[
c_t(V \otimes \C) = \sum_{k \geq 0} c_k(V \otimes \C) t^k  = c_t(V) \cdot c_t(\overline{V})
\]
equals 
\[
\prod (1 - x_i^2 t^2) = \sum (-1)^k e_k(x_i^2) t^{2k} \;,
\]
which expresses the Pontryagin classes
\[
p^\SO_k(V) = e_k(x_i^2)
\]
in terms of elementary symmetric functions of the Chern roots $x_i$ of $V \otimes \C$. \bigskip
 
\noindent
If $H_\G(z) := \hat{H}_\G(z^2)$ is an even power series, then the associated genus 
$\chi_\G$ of a $\C$-oriented manifold $M$ can be evaluated in terms of Pontryagin classes, since
\[
\prod \hat{H}_\G(x_i^2) := \BH_\G(p^\SO_k)
\]
for some polynomial $\BH_\G$; this factors $\chi_\G$ through a homomorphism
\[
\xymatrix{
M\U \ar[r] & M\SO \ar[r]^{\hat{\chi}_\G} & A[v] } \;.
\]
The complex vector bundle underlying a quaternionic vector bundle $V$, on the other hand,
can be decomposed as the sum of a complex bundle with its conjugate. In that case we have
\[
p^\SO_t(V) = p^\SO_t(W \oplus \overline{W}) = p^\SO_t(W)^2
\]
(at least, with coefficients in a $\Z[\half]$-algebra). The symplectic Pontryagin classes
of $V$ are defined by
\[
p^\Sp_t(V) = \sum (-1)^k c_{2k}(V) t^{2k}
\]
[20], so $p^\Sp_t(V) = p^\SO_t(W)$, hence $p^\SO_t(V) = (p^\Sp_t(W))^2$. Since $p^\SO_t(V)$
can be expressed in terms of the power sums $\sum x_i^{2k} = \ess^\SO_k$ of the Chern roots of
$V \otimes \C$ as 
\[
\exp(\; \sum \ess^\SO_{2k} \; \frac{t^{2k}}{k}) \;,
\]
we have 
\[
\ess^\SO_{2k} := \ess_{2k}(V \otimes \C) = 2 \ess_{2k}(V) := 2 \ess^\Sp_{2k}
\]
(in terms of the Chern roots of the complex structure underlying a quaternionic structure on 
$V$). \bigskip

\noindent
{\bf 2.3.1} Rewriting the logarithm of Weierstrass's product formula for $\Gamma$, we have
\[
\Gamma(1+z) = \exp(-\gamma z + \sum_{k > 1} \frac{\zeta(k)}{k} (-z)^k) \;;
\]
from this, and the duplication formula
\[
\Gamma(z) \Gamma(1-z) \; = \; \frac{\pi}{\sin \pi z}
\]
it follows that 
\[
\frac{x/2}{\sinh x/2} = \exp( \; \sum_{k \geq 1} \frac{\zeta(2k)}{(2\pi i)^{2k}} \; \frac{x^{2k}}{2k}) \;,
\]
with rational coefficients
\[
\frac{\zeta(2k)}{(2 \pi i)^{2k}} = - \frac{B_{2k}}{2(2k)!} \;.
\]
The $\hat{A}$-genus of an oriented manifold (corresponding to the group law in \S 1.3.2) can thus be
calculated by evaluating 
\[
\prod \Big( \frac{vx_i/2}{\sinh vx_i/2} \Big) = \exp( - \sum \frac{B_{2k}}{(2k)!} \; 
\frac{\ess^\SO_{2k}}{4k} \; v^{2k})
\]
on its fundamental class. If the manifold is $\bH$-oriented, this characteristic class equals the product
\[
\prod \Big( \frac{x_i/2}{\sinh x_i/2} \Big)^{1/2} 
\]
(now taken over the Chern roots of the complex bundle underlying the $\bH$-oriented structure). \bigskip

\noindent
{\bf Proposition.} {\it The genus of complex-oriented manifolds defined by the multiplicative series
\[
H_{\G_\infty}(x) = \Gamma(1+[\lambda]) = \Big(\frac{x/2}{\sinh x/2}\Big)^{1/2} \exp(i\frac{\gamma}{2\pi}x + 
\sum \frac{\zeta(\odd)}{(2 \pi i)^\odd} \; \frac{x^\odd}{\odd}) \in \C[[x]]
\]
agrees on the image of $M\Sp$ in $M\U$ with the $\hat{A}$-genus.} \bigskip

\noindent
[Because the odd terms in the exponential cancel, for a bundle of the form $W \oplus \overline{W}$.]
\bigskip

\noindent
{\bf 2.3.2} Note that the Witten genus [14]
\[
H_W(x) = \frac{x/2}{\sinh x/2} \prod_{n \geq 1} [(1 - q^ne^x)(1 - q^ne^{-x}]^{-1}
\]
can be written similarly, in terms of Eisenstein series, as 
\[
\exp(\sum_k G_{2k}(q) \; \frac{x^{2k}}{2k}) \;;
\]
but this deformation of the $\hat{A}$-genus is an {\bf even} function of $x$.
\bigskip 

\noindent
{\bf 2.4} The elementary symmetric functions $e_n$ and the corresponding power sums $\ess_n$ are
related by
\[
e(z) = \sum_{n \geq 0} e_n z^n := \prod_{k \geq 1}(1 + x_k z) = \exp(-\sum_{n \geq 1}
\frac{\ess_n}{n}(-z)^n) \;.
\]
The assignment $x_k \mapsto 1/k$ [5, 9, 14 I \S 2 ex 21] requires some care, but,
suitably interpreted, sends $\ess_k$ to $\zeta(k)$ if $k>1$, and $\ess_1$ to $\gamma$.
The formal power series
\[
{\rm Exp}_\infty(z) = z \cdot e(z)
\]
thus specializes to $\exp_\infty(z)$ under this mapping, defining a lift $\bG_\infty$ of $\G_\infty$
to a formal group law over the polynomial algebra $\Z[e_n \:| \: n \geq 1]$. Since its exponential
is defined over $\Z$, it is of additive type, and is in fact the universal such formal group law. \bigskip

\noindent
Similarly
\[
H_{\bG_\infty}(z) = \sum_{k \geq 0} h_k (-z)^k \;,
\]
in terms of the complete symmetric functions $h_k$. \bigskip

\section{The Real structure of $M\U$} \bigskip

\noindent
{\bf 3.1 Proposition.} {\it In the homotopy-commutative diagram
\[
\xymatrix{
M\U \ar[dr]^{\bgamma [\half]} \ar[r] & S^0[BU_+] \wedge \rH\Z \ar[r] & S^0[\Sp/\U_+ \wedge B\Sp_+] 
\wedge \rH\Z[\half] \ar[d]^{\zeta(\ev)} \\
S^0[\Sp/\U_+] \wedge M\Sp \ar[u] \ar[r] & S^0[\Sp/\U_+] \wedge \KO[\half] \ar[dr] \ar[r] & S^0[\Sp/\U_+] \wedge 
\rH \Q[v^{\pm 1}] \ar[d]^{\zeta(\odd)} \\
M\Sp \ar[u] \ar[r]^{\hat{A}} & \KO \ar[u] \ar[r] & \rH\C[v^{\pm 1}] }
\]
of spectra, the diagonal composition represents the $\Gamma$-genus.} \bigskip

\noindent
{\bf 3.2 Proof.} Here $S^0[G_+]$ is the suspension ring-spectrum defined by an $H$-space $G$, such as the 
fiber $\Sp/\U \; (\sim \Omega \Sp \sim B(\U/\Oh))$ of the quaternionification map $B\U \to B\Sp$. Note that 
the inclusion of the fiber into $B\U$ makes 
$S^0[B\U_+]$ (and hence $M\U$) into $S^0[\Sp/\U_+]$-modules). \bigskip

\noindent
The two vertical maps at the lower left side of the diagram are the obvious smash products with the unit
$S^0 \to S^0[\Sp/\U_+]$, while the horizontal maps across the middle of the diagram are smash products with 
the $\hat{A}$-genus, regarded as defined by the index of a Dirac operator on an $\bH$-oriented manifold,
followed by the Chern character on $\KO$. The top left-hand map is just the total characteristic number
homomorphisms of Boardman and Quillen, and can alternately be described as the composition
\[
M\U_* \to M\U_* \otimes S_* \to \Z \otimes S_* = S_*
\]
of the total Landweber-Novikov operation with Steenrod's cycle map
\[
1 \in H^0(B\U,\Z) \to H^0(M\U,\Z) = [M\U,\rH \Z]_0 \;.
\]
The (related) upper left-hand vertical and upper right-hand horizontal maps are more interesting. An element
of $M\Sp_*(\Sp/\U_+)$ can be interpreted as the bordism class of an $\bH$-oriented manifold $M$, equipped with a
map to $\Sp/\U$, and if we regard $M$ as merely complex-oriented, then the product composition
\[
M \to \Sp/\U_+ \wedge B\U_+ \to B\U_+
\]
defines a new complex orientation on $M$, and thus a ring homomorphism
\[
M\Sp_*(\Sp/\U_+) \to M\U_* \;.
\]
By [3], this is in fact an isomorphism away from the prime $(2)$; similarly, the composition
\[
\Sp/U_+ \wedge B\Sp_+ \to \Sp/\U_+ \wedge B\U_+ \to B\U_+
\]
defines an isomorphism 
\[
H_*(\Sp/\U,\Z[\half]) \otimes_{\Z[\half]} H_*(B\Sp,\Z[\half]) \cong H_*(B\U,\Z[\half])
\]
of Hopf algebras, which is the upper right-hand map. \bigskip
 
\noindent
Since the diagonal maps are defined by the diagram, only the right-hand vertical maps remain to be constructed,
but that is the content of \S 2.4: the power-sum generators of $H_*(B\U,\Q)$ map to normalized zeta-values
\[
\ess_k \mapsto \tilde{\zeta}(k) := (2 \pi i)^{-k} \zeta(k) \; {\rm if} \; k > 1 \; \; , \; \mapsto - \frac{\gamma}{2\pi} 
\cdot i \; {\rm if} \; k = 1 \;.
\]
This is factored into two steps:
\[
\zeta(\ev) : \ess_{2k} \mapsto \frac{B_{2k}}{4k(2k)!} \in \Q
\]
can be interpreted as defining the $\hat{A}$-genus, while
\[
\zeta(\odd) : \ess_{2k+1} \mapsto (-1)^{k+1} (2 \pi)^{-2k-1}\zeta(2k+1) \cdot i \;.
\]

\noindent
{\bf 3.3} Complex conjugation on $M\U$ is represented by the coordinate change $z \mapsto [-1](z)$ on the 
formal group, which corresponds to complex conjugation on the value group of the $\Gamma$-genus. In other 
words, the $\Gamma$-genus is naturally $\Z_2$-equivariant, with respect to the Galois action defined by 
the Real structure on complex cobordism. \bigskip

\noindent
Away from $(2)$, the Landweber-Novikov algebra of cobordism operations is an enveloping algebra of a $\Z_2$-graded
Lie (NB not super-Lie) algebra. The odd part corresponds, in classical Lie theory, to the tangent space
of the symmetric space associated to the complexification of a real Lie group; it acts transitively on
Spec $H_*(\Sp/\U,\Q)$, cf. [3, 17]. \bigskip

\section{Closing remarks} \bigskip

\noindent
{\bf 4.1} The index map $M\Sp \to \KO$ dates back to Conner and Floyd's 1968 work on the relation of cobordism
to $K$-theory, but seems to have received remarkably little attention: it is surely represented geometrically 
by a Dirac operator on $\bH$-oriented manifolds, but the question of a nice construction seems not to have
caught the differential geometers' attention. In view of this, I have not tried to define an explicit family 
of deformations of such an operator over $\Sp/\U$. \bigskip

\noindent
{\bf 4.2} R. Lu [8] has proposed an analytic interpretation of a variant of
the $\Gamma$-genus of a complex-oriented $M$ as a $\T$-equivariant Euler
class of its free loopspace, following Atiyah ([2]; see also [1]).
Lu's construction depends on a choice of {\bf polarization}
\[
\xymatrix{
\U/\Oh \ar[r] & {B\Gl_\res \sim B(\Z \times B\Oh)} \ar[d] \\
LM \ar@{.>}[u] \ar[ur] \ar[r] & {LB\U \sim B(L\U) \sim B(\Z \times B\U)}
}
\]
of the tangent bundle of $LM$: that is, a lift of the map classifying
its tangent bundle, to the restricted Grassmannian defined by writing
loops in the tangent space as a sum of something like positive and negative-frequency
components. Since $M$ is complex-oriented, such a lift exists, but is not
in general unique: it can be twisted by a map
\[
LM \to \U/\Oh \sim \Omega (\Sp/\U)
\]
[6 \S 2, 7, 10]. The free loops on a map $\alpha :M \to \Sp/\U \in \KO^2(X)$ thus define
a twist
\[
L(\alpha) : LM \to L(\Sp/\U) \sim \U/\Oh \times \Sp/\U \to \U/\Oh \;;
\]
its restriction to the subspace $M$ of constant loops defines a map to $\Sp/\U \sim \Omega \Sp$
which acts naturally on $\U/\Oh \sim \Omega (\Sp/\U)$, and it seems reasonable to expect that Lu's 
class for the polarized manifold $(M,L(\alpha))$ can be expressed in terms of $\bgamma(M)$ evaluated at
suitable values $\ess_{2k}(\alpha)$ of the deformation parameters. \bigskip

\noindent
{\bf 4.3} I have also not tried to pin down the two-local properties of $\bgamma$, which seem quite interesting. 
Away from $(2)$, $\Sp/\U$ is closely related [4] to $B\Sp(\Z)^+$, which is in turn related (via Siegel) to the 
$K$-theory spectrum of the symmetric monoidal category of Abelian varieties. This suggests that one might hope 
to see in the $\Gamma$-genus, some homotopy-theoretic residue of the intermediate Jacobians of complex projective 
manifolds. \bigskip

\noindent
{\bf 4.4} Kontsevich's original remarks were motivated by questions of quantization, and nothing
in the discussion above says much about that: homotopy theory is often revealing about the bones of a subject,
without resolving the surrounding analytical structures. \bigskip

\noindent
It is intriguing that the points $0, \; 1, \; \infty$ on the projective line seem to have naturally associated
genera and cohomology theories: the additive group at zero is related to de Rham theory, and the multiplicative
group at one to $K$-theory. The association of the point at infinity with the Kontsevich genus suggests it might
be related to a Galois theory of asymptotic expansions, along lines suggested by Cartier, Connes, Kreimer, Marcolli,
and others. \bigskip

\bibliographystyle{amsplain}

\end{document}